\documentclass[11pt,reqno]{amsart}

\usepackage[T1]{fontenc}
\usepackage[utf8]{inputenc}
\usepackage{mathpazo}
\usepackage[scaled=.95]{helvet}
\usepackage{microtype}
\usepackage[margin=1.10in]{geometry}
\usepackage{amsmath,amssymb,mathtools,mathrsfs}
\usepackage{xcolor}
\usepackage{enumitem}
\usepackage{hyperref}
\usepackage[nameinlink,noabbrev]{cleveref}
\usepackage[normalem]{ulem}

\hypersetup{
  colorlinks=true,
  linkcolor=blue!55!black,
  citecolor=blue!55!black,
  urlcolor=blue!55!black,
  pdfauthor={Juan J. Manfredi and Aelson Sobral},
  pdftitle={p-harmonic  polynomials}
}

\numberwithin{equation}{section}
\theoremstyle{plain}
\newtheorem{theorem}{Theorem}[section]

\newtheorem{lemma}[theorem]{Lemma}

\theoremstyle{definition}

\theoremstyle{remark}
\newtheorem{remark}[theorem]{Remark}

\newcommand{\Rn}{\mathbb R^n}
\newcommand{\Sn}{\mathbb S^{n-1}}
\newcommand{\norm}[1]{\left\lVert #1\right\rVert}
\newcommand{\abs}[1]{\left\lvert #1\right\rvert}
\newcommand{\diver}{\operatorname{div}}
\newcommand{\E}{\mathcal E}

\setlist[enumerate]{leftmargin=2.2em,itemsep=4pt,topsep=4pt}
\setlist[itemize]{leftmargin=2em,itemsep=3pt,topsep=4pt}
\allowdisplaybreaks[2]

\title[$p$-harmonic  polynomials]{$p$-harmonic  polynomials}

\author[J. J. Manfredi]{ Juan J. Manfredi}
\address{Department of Mathematics, University of Pittsburgh,
301 Thackeray Hall, Pittsburgh, PA 15260, USA}
\email{manfredi@pitt.edu}

\author[A. Sobral]{Aelson Sobral}
\address{Applied Mathematics and Computational Sciences (AMCS),
Computer, Electrical and Mathematical Sciences and Engineering Division
(CEMSE), King Abdullah University of Science and Technology (KAUST),
Thuwal, 23955-6900, Kingdom of Saudi Arabia}
\email{aelson.sobral@kaust.edu.sa}

\date{}

\begin{document}

\begin{abstract}
We prove that every real polynomial $p$-harmonic function on $\Rn$ is affine when $n\ge3$ and $p\ge n-1/13$. The proof uses the Bernstein--Szeg\H{o} inequality on great circles and a deficit estimate at the first Chebyshev node.
\end{abstract}

\maketitle

\section{Introduction and statement of the result}

Let $1<p<\infty$, $p\ne2$, and consider weak solutions of the $p$-Laplace equation
\[
\Delta_pu=\diver\bigl(\abs{Du}^{p-2}Du\bigr)=0 \qquad\text{in }\Rn.
\]
Whether this equation admits a homogeneous polynomial solution of degree at least two is a longstanding open question, raised by Lewis \cite{Lewis1980} and recorded in its present form by Tkachev \cite{TkachevScilag}. Lewis \cite{Lewis1980} ruled out homogeneous polynomial solutions of degree at least two in the plane, as well as solutions of degree two in every dimension. Tkachev \cite{Tkachev} treated degree three by algebraic methods, and Lewis and Vogel \cite{LewisVogel} treated degree four, as well as degree five in three variables. In this paper we give an explicit range of $p$, depending only on the dimension, in which nonexistence holds in every degree.

Write $\Sn=\{x\in\Rn:\abs{x}=1\}$ and set
\begin{equation}\label{eq:E}
\E_p(u)=\abs{Du}^2\Delta u+(p-2)\langle D^2u\,Du,Du\rangle.
\end{equation}
If a nonconstant polynomial $u$ is a weak solution of $\Delta_pu=0$, then $u$ is a classical solution on the open set $\{Du\ne0\}$, so that $\E_p(u)$ vanishes there. At points where $Du=0$, both terms in \eqref{eq:E} vanish. Thus $\E_p(u)=0$ throughout $\Rn$.

\begin{theorem}\label{thm:main}
Let $n\ge3$ and $p\ge n-1/13$. Every real polynomial $u$ on $\Rn$ satisfying $\E_p(u)=0$ is affine. In particular, every polynomial $p$-harmonic function on $\Rn$ is affine.
\end{theorem}

\begin{remark}\label{rem:ptwo}
The theorem rules out non-affine polynomial $p$-harmonic functions throughout the range $p\ge n$, with no restriction on the degree. This gives a finite-$p$ counterpart of the nonexistence of $\infty$-harmonic homogeneous polynomials of degree at least two. See, for example \cite{YifengYu}. The deficit estimate also allows the argument to cover the interval $n-1/13\le p<n$. The hypotheses ensure $p>2$, as required in the proof.
\end{remark}

\begin{remark}\label{rem:constant}
The constant $1/13$ is chosen to keep the estimates explicit and the proof simple. A sharper evaluation of the final numerical estimates gives only a small improvement, with a threshold of approximately $0.0789$. We do not claim that $1/13$ is optimal for the problem.
\end{remark}

The plan of the paper is as follows. In \Cref{sec:deficit} we prove a deficit estimate at the first Chebyshev node (\Cref{lem:deficit}). In \Cref{sec:proof}, we prove \Cref{thm:main} by evaluating the equation at a well-chosen point on the sphere, $\theta=\pi/m$.

\section{A deficit estimate}\label{sec:deficit}

The Bernstein--Szeg\H{o} inequality states that a real trigonometric polynomial $T$ of degree at most $m\ge1$, with $\norm{T}_\infty\le1$, satisfies $T^2+(T'/m)^2\le1$; it goes back to Szeg\H{o} and to van der Corput and Schaake \cite{vanderCorputSchaake}. A proof can be found in \cite[Theorem~3.7]{QueffelecZarouf}. Applying it also to $T'/m$ gives
\begin{equation}\label{eq:bernstein}
T^2+\frac{(T')^2}{m^2}\le1, \qquad \frac{(T')^2}{m^2}+\frac{(T'')^2}{m^4}\le1.
\end{equation}
In particular, $\norm{T''}_\infty\le m^2$.

The proof of the following lemma is an interpolation argument at the Chebyshev nodes, in the spirit of the proof by Duffin and Schaeffer of the Markov brothers' inequality \cite{DuffinSchaeffer}.

\begin{lemma}\label[lemma]{lem:deficit}
Let $m\ge2$ be an integer, and let $T$ be a real trigonometric polynomial of degree at most $m$, with $\norm{T}_\infty\le1$ and $T(0)=1$. If $d=m^2+T''(0)$, then
\begin{equation}\label{eq:deficit-bound}
0\le1+T(\pm\pi/m)\le d\bigl(1-\cos(\pi/m)\bigr)
\le\frac{\pi^2d}{2m^2}.
\end{equation}
\end{lemma}

\begin{proof}
Write $P(\cos t)=(T(t)+T(-t))/2$, where $P$ is a polynomial of degree at most $m$, and let $C_m(\cos t)=\cos(mt)$. Then $\abs{P}\le1$ on $[-1,1]$ and $P(1)=1$. Differentiating $P(\cos t)$ twice at $t=0$ gives $P'(1)=-T''(0)$, and $C_m'(1)=m^2$, so that $d=C_m'(1)-P'(1)$. At the nodes $x_j=\cos(j\pi/m)$, $0\le j\le m$, differentiation of the Lagrange interpolant of $C_m-P$ gives
\[
d=\sum_{j=1}^m a_j\bigl(1-(-1)^jP(x_j)\bigr),
\qquad
a_j=\begin{cases}
2/(1-x_j),&1\le j<m,\\
1/2,&j=m.
\end{cases}
\]
The $j=0$ term vanishes because $(C_m-P)(1)=0$. Indeed, for $w(x)=(x^2-1)C_m'(x)$, the Lagrange cardinal polynomial $\ell_j$ satisfies
\[
\ell_j'(1)=\frac{w'(1)}{w'(x_j)(1-x_j)}=(-1)^ja_j,
\qquad 1\le j\le m.
\]
The last equality follows by differentiating $w(\cos t)=-m\sin t\sin(mt)$, taking limits at the endpoints. All summands are nonnegative. Taking $j=1$ and writing $h=\pi/m$, we get
\[
d\ge\frac{(1+T(h))+(1+T(-h))}{1-\cos h}.
\]
Each term in the numerator is nonnegative. The conclusion follows from $1-\cos h\le h^2/2$.
\end{proof}

\section{Proof of the theorem}\label{sec:proof}

\begin{proof}[Proof of \Cref{thm:main}]
Suppose that $u$ has degree $m\ge2$. The component of degree $3m-4$ in $\E_p(u)$ is $\E_p(u_m)$, where $u_m$ is the highest homogeneous part of $u$. We may therefore assume that $u$ is homogeneous of degree $m$ and, after replacing $u$ by $-u$ if necessary, normalize it so that
\[
u(e)=\max_{\Sn}\abs{u}=1.
\]
For a unit vector $v\perp e$, put $T_v(t)=u(e\cos t+v\sin t)$. This is a trigonometric polynomial of degree at most $m$. At $e$, homogeneity and maximality give $Du(e)=me$ and $D^2u(e)[e,e]=m(m-1)$. Since $T_v''(0)=D^2u(e)[v,v]-\langle Du(e),e\rangle$, we obtain
\[
d_v\coloneqq m^2+T_v''(0)
=D^2u(e)[v,v]+m(m-1)\ge0,
\]
by \eqref{eq:bernstein}. Evaluating $\E_p(u)=0$ at $e$ and summing over any orthonormal basis $v_1,\ldots,v_{n-1}$ of $e^\perp$ gives
\begin{equation}\label{eq:trace}
\sum_{i=1}^{n-1}d_{v_i}=(n-p)m(m-1).
\end{equation}
Since every $d_{v_i}$ is nonnegative, $p\le n$, and the theorem is proved for $p>n$. With $\delta=n-p$, we have $0\le\delta\le1/13$ and $d_v\le\delta m(m-1)$ for every $v$.

Write $F(\theta,v)=u(e\cos\theta+v\sin\theta)$ and $\theta_0=\pi/m$. By \Cref{lem:deficit},
\begin{equation}\label{eq:latitude}
-1\le F(\theta_0,v)\le-1+\frac{\pi^2\delta}{2}
\qquad\text{for every }v.
\end{equation}
Choose $v_0$ maximizing $F(\theta_0,\cdot)$ on the unit sphere in $e^\perp$. At this point $D_vF=0$ and $\Delta_vF\le0$, where these operators act on that sphere. Set
\[
a=-F(\theta_0,v_0),\qquad b=\frac{F_\theta}{m},\qquad
c=\frac{F_{\theta\theta}}{m^2},\qquad q=a^2+b^2.
\]
All derivatives below are evaluated at $(\theta_0,v_0)$. Since $\pi^2<10$, \eqref{eq:bernstein} and \eqref{eq:latitude} imply
\[
1\ge a\ge1-\frac{\pi^2\delta}{2}>\frac8{13},\qquad
q\le1,\qquad b^2+c^2\le1.
\]
Put $s=\sqrt{1-b^2}$, so that $c\le s$ and $s\ge a>0$.

In the coordinates $x=r(e\cos\theta+v\sin\theta)$, with $r>0$, $0<\theta<\pi$, and $v$ in the unit sphere of $e^\perp$, the Laplacian reads
\begin{equation}\label{eq:suspension}
\Delta u=u_{rr}+\frac{n-1}{r}\,u_r
+\frac1{r^2}\Bigl(u_{\theta\theta}+(n-2)\cot\theta\,u_\theta
+\frac{\Delta_vu}{\sin^2\theta}\Bigr).
\end{equation}
Since $D_vF=0$, the gradient at $x_0=e\cos\theta_0+v_0\sin\theta_0$ has only radial and polar components, and $\abs{Du}^2=m^2q>0$; in particular, $Du(x_0)$ lies in the plane $\Pi=\operatorname{span}\{e,v_0\}$, so that only the restriction of $D^2u$ to $\Pi$ enters $\langle D^2u\,Du,Du\rangle$. Because $\Pi$ is an affine plane, the second derivatives of $u$ along $\Pi$ are those of the restriction $r^mF(\theta,v_0)$, computed in polar coordinates: in the orthonormal frame $\{e_r,e_\theta\}$ at $x_0$,
\begin{equation}\label{eq:frame}
Du=mF\,e_r+F_\theta\,e_\theta,
\qquad
D^2u\big|_\Pi=
\begin{pmatrix}
m(m-1)F & (m-1)F_\theta\\[2pt]
(m-1)F_\theta & F_{\theta\theta}+mF
\end{pmatrix}.
\end{equation}
Substituting $F=-a$, $F_\theta=mb$, $F_{\theta\theta}=m^2c$ in \eqref{eq:suspension} and \eqref{eq:frame}, and dividing $\E_p(u)=0$ by $m^2\abs{Du}^2=m^4q$, we obtain
\begin{equation}\label{eq:polar}
\begin{aligned}
0=\frac{\E_p(u)}{m^2\abs{Du}^2}
={}&c-\left(p-1+\frac\delta m\right)a
 +\frac{(n-2)\cot\theta_0}{m}b\\
&+(p-2)\frac{b^2(c-a)}q
 +\frac{\Delta_vF}{m^2\sin^2\theta_0}.
\end{aligned}
\end{equation}
Since $p>2$, the coefficient of $c$ is positive. Using $c\le s$, $\Delta_vF\le0$, and $0\le\cot\theta_0/m\le1/\pi$, we obtain
\[
0\le s-(p-1)a+\frac{(n-2)\abs b}{\pi}
 +(p-2)\frac{b^2(s-a)}q.
\]
Moreover, $s^2-a^2=1-q$ and $q+a^2/q\ge2a$ give
\[
0\le\frac{b^2(s-a)}q
=\frac{1+a^2-q-a^2/q}{s+a}
\le\frac{(1-a)^2}{2a}\eqqcolon A.
\]
Put $k=n-2\ge1$. Since $p-1=k+1-\delta$, $p-2=k-\delta$, and $s-(1-\delta)a\ge0$, it follows that
\[
\begin{aligned}
0&\le s-(p-1)a+k\frac{\abs b}{\pi}+(p-2)A\\
 &\le k\left(s+\frac{\abs b}{\pi}-(2-\delta)a+A\right).
\end{aligned}
\]
Finally, Cauchy--Schwarz and $\pi>3$ give $s+\abs b/\pi\le\sqrt{1+\pi^{-2}}<17/16$. Using $\delta\le1/13$ and $a>8/13$, we conclude that
\[
0\le s+\frac{\abs b}{\pi}-(2-\delta)a+A
<\frac{17}{16}-\frac{200}{169}+\frac{25}{208}
=-\frac1{1352}<0,
\]
a contradiction.
\end{proof}
\section{Declaration of generative AI and AI-assisted technologies in the manuscript preparation process} During the preparation of this work, the author(s) used LLM models from OpenAI and Anthropic for assistance in the discovery and verification of proofs. The author(s) reviewed and edited the output as needed and take full responsibility for the content of the published article. All mathematical results presented in this paper have been formally verified in Lean~4 using the Mathlib library. The Lean code and an explanation of the formalization process are available at \url{https://github.com/na0b/pharmonicpoly-lean}.

\bigskip

{\small \noindent{\bf Acknowledgments.} JM was partially supported by a gift from the Simons Foundation (Gift Number 962828). This publication is based upon work supported by King Abdullah University of Science and Technology (KAUST). AS acknowledges the hospitality of the mathematics department at the University of Pittsburg, where part of this research was conducted.}

\end{document}